\documentclass{amsart}
\usepackage[dvips]{graphicx}
\usepackage{amscd}
\usepackage{amsmath}
\usepackage{amsxtra}
\usepackage{amsfonts}
\usepackage{amssymb}
\newtheorem{theorem}{Theorem}[section]

\newtheorem{proposition}[theorem]{Proposition}
\theoremstyle{definition}
\newtheorem{definition}[theorem]{Definition}
\newtheorem{remark}[theorem]{Remark}

\theoremstyle{remark}

\renewcommand{\theclaim}{\textup{\theclaim}}

\newtheorem*{acknowledgements}{Acknowledgements}

\numberwithin{equation}{section}

\def\openone
{\mathchoice
{\hbox{\upshape \small1\kern-3.3pt\normalsize1}}
{\hbox{\upshape \small1\kern-3.3pt\normalsize1}}
{\hbox{\upshape \tiny1\kern-2.3pt\SMALL1}}
{\hbox{\upshape \Tiny1\kern-2pt\tiny1}}}
\makeatletter
\newbox\ipbox
\newcommand{\ip}[2]{\left\langle #1\mathrel{\mathchoice
{\setbox\ipbox=\hbox{$\displaystyle \left\langle\mathstrut #1#2\right\rangle$}
\vrule height\ht\ipbox width0.25pt depth\dp\ipbox}
{\setbox\ipbox=\hbox{$\textstyle \left\langle\mathstrut #1#2\right\rangle$}
\vrule height\ht\ipbox width0.25pt depth\dp\ipbox}
{\setbox\ipbox=\hbox{$\scriptstyle \left\langle\mathstrut #1#2\right\rangle$}
\vrule height\ht\ipbox width0.25pt depth\dp\ipbox}
{\setbox\ipbox=\hbox{$\scriptscriptstyle \left\langle\mathstrut #1#2\right\rangle$}
\vrule height\ht\ipbox width0.25pt depth\dp\ipbox}
} #2\right\rangle}
\newcommand{\diracb}[1]{\left\langle #1\mathrel{\mathchoice
{\setbox\ipbox=\hbox{$\displaystyle \left\langle\mathstrut #1\right.$}
\vrule height\ht\ipbox width0.25pt depth\dp\ipbox}
{\setbox\ipbox=\hbox{$\textstyle \left\langle\mathstrut #1\right.$}
\vrule height\ht\ipbox width0.25pt depth\dp\ipbox}
{\setbox\ipbox=\hbox{$\scriptstyle \left\langle\mathstrut #1\right.$}
\vrule height\ht\ipbox width0.25pt depth\dp\ipbox}
{\setbox\ipbox=\hbox{$\scriptscriptstyle \left\langle\mathstrut #1\right.$}
\vrule height\ht\ipbox width0.25pt depth\dp\ipbox}
}\right. }
\newcommand{\dirack}[1]{\left. \mathrel{\mathchoice
{\setbox\ipbox=\hbox{$\displaystyle \left.\mathstrut #1\right\rangle$}
\vrule height\ht\ipbox width0.25pt depth\dp\ipbox}
{\setbox\ipbox=\hbox{$\textstyle \left.\mathstrut #1\right\rangle$}
\vrule height\ht\ipbox width0.25pt depth\dp\ipbox}
{\setbox\ipbox=\hbox{$\scriptstyle \left.\mathstrut #1\right\rangle$}
\vrule height\ht\ipbox width0.25pt depth\dp\ipbox}
{\setbox\ipbox=\hbox{$\scriptscriptstyle \left.\mathstrut #1\right\rangle$}
\vrule height\ht\ipbox width0.25pt depth\dp\ipbox}
} #1\right\rangle}

\newcommand{\linfr}{L^{\infty}\left(\mathbb{R}\right)}

\input cyracc.def

\begin{document}
\title[The Canonical Equivalence Relation of the Thompson Group]{The von Neumann Algebra of the Canonical Equivalence Relation of the Generalized Thompson Group }
\author{Dorin Ervin Dutkay}
\address{Department of Mathematics\\
The University of Iowa\\
14 MacLean Hall\\
Iowa City, IA 52242-1419\\
U.S.A.}
\email{ddutkay@math.uiowa.edu}
\author{Gabriel Picioroaga}
\address{Department of Mathematics\\
The University of Iowa\\
14 MacLean Hall\\
Iowa City, IA 52242-1419\\
U.S.A.}
\email{gpicioro@math.uiowa.edu}
\thanks{}

\begin{abstract} We study the equivalence relation $R_N$ generated by the (non-free) action of 
the generalized Thompson group $F_N$ on the unit interval. We show that this relation is a standard,
quasipreserving ergodic equivalence relation. Using results of 
Feldman-Moore, Krieger and Connes we prove that the von Neumann
algebra $M(R_N)$ associated to $R_N$ is the hyperfinite type
$III_{\lambda}$ factor, with $\lambda=1/N$. 
\\
Moreover we analyze $R_N$ and $F(N)$ in connection with Gaboriau's work on costs
of groups. We prove that the cost $C(F(N))=1$ for any $N\geq 2$ and for $N=2$ we precisely find 
a treeing of $R_N$.
\end{abstract}\maketitle

\section{\label{Intro}Introduction}
In the following we prepare the definitions we need in this paper. We also mention some known results we are going to use: we
follow \cite{Gab}, \cite{FM}, \cite{Br} and \cite{Can}.
\\ We say that  $R$ is a SP1 equivalence relation on a standard probability space $(X,\lambda)$ if
\par (S) Almost each orbit $R[x]$ is at most countable and $R$ is a Borel subset of $X\times X$.
\par (P) For any $T\in\mbox{Aut}(X,\lambda)$ such that $\mbox{graph}T\subset R$ we have that $T$ preserves the measure $\lambda$.
\\ We say that $R$ is standard if only (S) is satisfied. Also, $R$ is called quasi-preserving if the saturation
(through $R$) of a null set is null. \\
From now on, unless specified otherwise, each equivalence relation satisfies (S).
Next we define  "graphing" and "treeability" with respect to $R$.
This is just a simple adaptation of the SP1 situation (see \cite{Gab}).
\begin{definition}\label{def1}i) A countable family $\Phi=(\varphi_i:A_i\rightarrow B_i)_{i\in I}$
of Borel partial isomorphisms between Borel subsets of $(X,\lambda)$ is called a
graphing on $(X,\lambda)$ (we do not require that the $\varphi_i$'s preserve $\lambda$).\\
ii) The equivalence relation $R_{\Phi}$ generated by a graphing $\Phi$ is the smallest equivalence relation $S$ such that
$(x,y)\in S$ iff $x$ is in some $A_i$ and $\varphi_i(x)=y$. \\
iii) An equivalence relation $R$ is called treeable if there is a graphing $\Phi$ such that $R=R_{\Phi}$ and almost
every orbit $R_{\Phi}[x]$ has a tree structure. In such case $\Phi$ is called a treeing of $R$.\\
iv) $R$ is ergodic iff any saturated Borel set has measure $0$ or $1$.
\end{definition}
\begin{remark} For (SP1) $R$'s the same notions are considered in \cite{Gab} provided the $\varphi_i$'s preserve the measure.
One can consider the quantity $C(\Phi)=\sum\lambda(A_i)$. The cost
of a (SP1) equivalence relation will simply be
$$C(R)=\inf\{C(\Phi) |  \Phi\mbox{ is a graphing of } R\}.$$ It is
the preserving property that allows one to conclude the infimum is
attained iff $R$ admits a treeing (see Prop.I.11 and Thm.IV.1 in
\cite{Gab}). Next
Gaboriau defines the cost of a discrete countable group $G$ as \\
$$C(G)=\inf\{C(R) |  R \mbox{  coming from a free, preserving
action of }G\mbox{ on }X \}.$$  We highlight the following result that 
"measures" the non-amenability of cost 1 groups (any amenable
group has cost=1):
\par $Theorem$(\cite{Gab}, Corollaire VI.22) $Any$ $non-amenable$ $cost$ $1$
$group$ $is$ $anti-treeable$ (i.e., any SP1 equivalence relation
coming from a free action
is not treeable).\\

 Among many examples of groups whose costs are calculated, the Thompson group is
shown to have cost=1 (using the infinite presentation of the group
and one of the tools developed by Gaboriau). Now any countable
discrete group comes with a free preserving action on some
standard probability space, namely the Bernoulli shifts (thus the
infimum in $C(G)$ does make sense). However, to handle (in terms
of (non)treeability) the SP1 relation determined by this purely
theoretical action may be very hard. Until we find a suitable
action of the generalized Thompson group, we are content to study its 
canonical action on $([0,1],\lambda)$, where $\lambda$ denotes the Lebesgue measure. 
Certainly this is not a (SP1) relation but it is (S) and quasi-preserving. 
\end{remark}
Let us  introduce some basics facts about the Thompson groups.
\begin{definition}\label{t}The Thompson group $F$ is the set of piecewise linear homeomorphisms from the closed unit interval $[0,1]$
to itself that are differentiable except at finitely many dyadic rationals and such that on intervals of differentiability
the derivatives are powers of $2$. \\
If $N\geq 2$, replacing above the dyadic rationals by $N$-adic rationals and the power of 2 slopes by 
powers of $N$, we obtain one of the generalized versions of the Thompson group. We will denote it by $F(N)$. 
\end{definition}
\begin{remark} It is shown that $F(N)$ above is a countable subgroup of the group of all homeomorphisms from $[0,1]$ to $[0,1]$.
Two presentations of $F$ are found. One finite presentation comes from the fact that $F$ is generated by the functions $A$
and $B$ defined below
$$A(x)=\left\{
\begin{array}{lr}
x/2,&\mbox{ } 0\leq x\leq 1/2\\
x-1/4,&\mbox{ } 1/2\leq x\leq 3/4\\
2x-1,&\mbox{ } 3/4\leq x\leq 1
\end{array}\right.\mbox{   ,    }
B(x)=\left\{\begin{array}{lr}
x,&\mbox{ }0\leq x\leq1/2\\
x/2+1/4,&\mbox{ }1/2\leq x\leq 3/4\\
x-1/8,&\mbox{ }3/4\leq x\leq 7/8\\
2x-1,&\mbox{ }7/8\leq x\leq 1
\end{array}\right. $$
The relations between generators $A$ and $B$ are $[AB^{-1},A^{-1}BA]=1$ and \\
$[AB^{-1},A^{-2}BA^2]=1]$. $F(N)$ has also a finite presentation, see \cite{Br}. However, 
for computing the cost of $F(N)$ we will make use  of the following infinite presentation 
 $$F(N)=\left< x_0, x_1, ...x_i,...|\mbox{ }x_jx_i=x_ix_{j+N-1}\mbox{, }i<j\mbox{ }\right>$$

\end{remark}
Next we will introduce the von Neumann algebra of an equivalence
relation. We follow \cite{FM} in the particular case when the
2-cocycle $\sigma$ is trivial. Let $R$ be a standard equivalence
relation on the standard probability space $(X,\lambda)$. The
Hilbert space the algebra acts upon is $H:=L^2(R,\nu_r)$ where
$\nu_r$ is the right-counting measure on $R$. When no confusion
with the left-counting measure $\nu_l$ may arise, we will write
 $\nu$ instead of $\nu_r$.  E.g., if $f\in H$, its squared norm is given by
$$\int |f(x,y|^2d\nu(x,y)=\int (\sum_{(z, x)\in R}|f(x,z)|^2)d\lambda(x)$$
For a left-finite function $a:R\rightarrow\mathbb{C}$, we denote
by $L_a$ the bounded operator
$$L_a\varphi(x,y)=\sum_za(x,z)\varphi(z,y).$$
Then $M(R)$ is defined as $\{L_a | a\mbox{   is left-finite
}\}^{''}$. It is known that $L^{\infty}(X)$ can be embedded as a
Cartan subalgebra into $M(R)$. Also, $\varphi_0$ the
characteristic function of the diagonal in $R$ is a separating and
cyclic vector. Any element $L\in M(R)$ can be written as
$L_{\psi}$ (where $\psi=L\varphi_0$), meaning that
$$L\varphi(x,y)=\sum_z\psi (x,z)\varphi(z,y)$$ for all $\varphi\in
H$ and all $(x,y)\in R$. Now the multiplication on $M(R)$ can be
written as a convolution over $R$:
$L_{\psi_1}*L_{\psi_2}=L_{\psi_1*\psi_2}$ where
$$\psi_1*\psi_2(x,y)=\sum_z\psi_1(x,z)\psi_2(z,y),\quad(x,y)\in R.$$
Moreover, if $R$ is ergodic then $M(R)$ is a factor. \\
If the equivalence relation $R$ is coming from 
the action of a discrete countable group  $G$ on the probability space $(X,\lambda)$, then $M(R)$ is 
the crossed product of $L^{\infty}(X)$ by $G$. This is exactly the situation we will work in, however 
we prefer the Feldman-Moore setting.   
\par
It is easy to show that if the measure is $R$-invariant (i.e., $R$
satisfies (P) above) then the state $<\mbox{ }_{\mbox{
}}^{.}\varphi_0, \varphi_0>$ is a trace; in this case $M(R)$ is a
factor of type II. If there is no $\sigma$-finite measure $\mu$,
$R$-invariant such that $\mu\prec\lambda$ then by Theorem 2.4 in
\cite{Kr}, $M(R)$ has to be of type III. (The non-existence of
such $\mu$ proves that $R$ is of type III, see the terminology in
\cite{FM1}). This is the result we are going to use for the
canonical $R_N$ on $F(N)$. Also, it turns out that the Connes spectrum
has a nice description for factors coming from ergodic equivalence
relations, namely the asymptotic range of the map $D:R\rightarrow
\mathbb{R}_{+}$
 $$D(x,y)=\frac{\partial\nu_l}{\partial\nu_r}(x,y)$$ therefore, to get the type of the factor suffices 
to "compute" the values of $D$. Following  \cite{Co} we obtain that $M(R_N)$ is the crossed-product 
of the hyperfinite $II_{\infty}$ factor by $\mathbb{Z}$.
\section{\label{m(r)}$M(R_N)$}
In the following if a "measurable" statement is made with respect to points on the real line (or plane) then it is 
understood that the measure taken into account is the Lebesgue measure. 

\begin{definition}The equivalence relation $R_N\subset [0,1]\times [0,1]$ defined by $(x,y)\in R_N$ iff there exists $f\in F(N)$ 
such that $f(x)=y$ is called the canonical equivalence relation of 
the generalized Thompson group. 
\end{definition} 
\begin{remark}\label{r1}We can work with a subrelation of $R_N$ (still denoted $R_N$) which is $R_N$ except a set of (product) 
measure $0$. This change will not affect the (S) or (P) properties nor the construction of $M(R_N)$. In our case $R_N$ is
replaced by $R_N$ minus the points of rational coordinates. 
\end{remark}
\begin{remark}\label{r2} We pause for a moment to distinguish the situations $N=2$ and $N>2$. Let $R$ be the equivalence relation 
generated by the $ax+b$ group with $a$ of the form $2^n$,
$n\in\mathbb{Z}$ and $b\in\mathbb{R}$, dyadic, i.e. $b=k/2^m$ for
some $k,m\in\mathbb{Z}$. Notice that $R_2$ is the restriction of $R$ to the unit square (however, this is not obvious: given 
$y=ax+b$ one has to construct $f\in F$ such that $f(x)=y$ and this can be carried-out by using the properties 
of the Thompson group, see Lemma 4.2 in \cite{Can}). Interestingly enough, a similar "localization" for $R_N$ with odd 
$N$ fails to hold true: we prove that for any $x\in[0,1]\setminus\mathbb{Q}$ there exists no $f\in F(N)$ such that  
$f(x)=x+\frac{k}{N^p}$ where $k$ is odd and $p\geq0$.
\par
First, take $x$ $N$-adic, $x=\frac{a}{N^r}$  with $r\geq p$. Then
$$x+\frac{k}{N^p}=\frac{a+kN^{r-p}}{N^r}=:\frac{b}{N^r},$$
and observe that $a$ and $b$ have different parity.
\par
Assume now that there is $f\in F(N)$ such that
$f(x)=x+\frac{k}{N^p}$. Consider all the points of
non-differentiability $\{l_i/N^{s_i}\,|\, i\in\mathcal{I}\}$ of
$f$ and take $r\geq\max\{s_i\,|\, i\in\mathcal{I}\}$.
\par
Let $N^{q_i}$ ($q_i\in\mathbb{Z}$) be the slope of $f$ on the
interval $[i/N^r,(i+1)/N^r]$, \\ $(i\in\{0,...,N^r-1\})$. Then
$$f\left(\frac{i+1}{N^r}\right)=f\left(\frac{i}{N^r}\right)+N^{q_i}\frac{1}{N^r},\quad(i\in\{0,...,N^r-1\}).$$
By induction, since $f(0)=0$ we obtain that
$$f\left(\frac{k}{N^r}\right)=\sum_{i=0}^{k-1}N^{q_i}\frac{1}{N^r},\quad(k\in\{1,...,N^r\}).$$
In particular

$$\frac{b}{N^r}=f\left(\frac{a}{N^r}\right)=\sum_{i=0}^{a-1}N^{q_i}\frac{1}{N^r}$$

Now take $q\geq\max\{-q_i\,|\,i\in\{1,...,a-1\}\}$ and multiply by
$N^q$:
\begin{equation}\label{eqrem1}
bN^q=\sum_{i=0}^{a-1}N^{q+q_i}.
\end{equation}
Since $a$ and $b$ have different parity and since $N$ is odd, it
follows that the terms of the equality (\ref{eqrem1}) have
different parity. This is a contradiction which shows that $d$ is
not equivalent to $d+k/N^p$ for $d$ $N$-adic.
\par
Now, if $x$ is not rational, assume that $x$ is equivalent to
$x+k/N^p$. This means that there is an $f\in F(N)$ such that
$f(x)=x+k/N^p$. Take a small interval around $x$ where $f$ is
differentiable. On this interval $f$ has the form $f(y)=N^sy+e$
with $s\in\mathbb{Z}$ and $e$ $N$-adic. But then $s=0$ and
$e=k/N^r$ otherwise $x+k/N^r=N^sx+e$ and this would imply that $x$
is rational. So on this interval $f(y)=y+k/N^p$. We can find an
$N$-adic point in this interval, call it $d$, such that $f(d)=d+k/N^p$
and this contradicts the fact that $d$ and $d+k/N^p$ are not
equivalent. 
\end{remark}
 
We have to make sure that $R_N$ is standard. The finite presentation of $F(N)$ implies that it is 
quasi-preserving. It is not hard to see 
that the group $F(N)$ is at most countable: given $x_1$, $x_2$..., $x_k$ a list of N-adic points in 
$[0,1]$ and a list of power of $N$ slopes there can be at most one element $f\in F(N)$ that fullfils these 
data. Therefore $F(N)$ is at most countable. We will actually show it is countable by displaying a 
non-trivial element in $F(N)$, useful also in the proofs below.
\begin{proposition}Let $d$ a N-adic in $[0,1]$ and $p\in\mathbb{Z}$ such that $d<N^p$. Then the following 
function is an infinite order element of $F(N)$:
$$A_{d,p}(x)=\left\{
\begin{array}{lr}
x/N^p,&\mbox{ } 0\leq x\leq d\\
x-d+d/N^p,&\mbox{ } d\leq x\leq 1-d/N^p\\
N^px+1-N^p,&\mbox{ } 1-d/N^p\leq x\leq 1
\end{array}\right. $$
\end{proposition}     
\begin{proof} The way $A_{d,p}$ is defined shows that it is an element of $F(N)$. Also, $A_{d,p}\neq$id, 
therefore all its iterates are distinct elements of $F(N)$. 
\end{proof} 
 We will show that the von Neumann algebra $M(R_N)$ is the type $III_{1/N}$ 
hyperfinite factor. We first prove ergodicity in order to insure that we are dealing with a factor.

\begin{proposition}\label{p1} The equivalence relation $S_N$ defined on [0,1] by $(x,y)\in S_N$ iff 
there exists $f\in F(N)$ such that $f(x)=y$ and $f^{'}(x)=1$, is an ergodic subrelation of $R_N$. Moreover, 
$S_N$ is a (SP1) hyperfinite equivalence relation with infinite orbits.
\end{proposition}
\begin{proof} Notice that if $(x,y)\in S_N$ through some $f\in F(N)$ then $f^{'}=1$ on a neighborhood of 
$x$ ($x$ not being N-adic). Clearly $S_N\subset R_N$. Let now $X$ be a $S_N$-saturated set. We show that 
for any $0<d_1<d_2<1$ N-adic numbers the following equality holds:
\begin{equation}\label{e4}
\lambda(X\cap [d_1,d_2])=\lambda(X\cap [0,d_2-d_1])
\end{equation}
Choose $p\in\mathbb{N}$ large enough such that $d_2<1-d_1/N^p$. Because\\
 $[d_1,d_2]\subset [d_1,1-d_1/N^p]$ and $A_{d_1,p}$ has slope 1 on $[d_1,1-d_1/N^p]$ we have 
$$\lambda(X\cap [d_1,d_2])=\lambda(A_{d_1,p}(X\cap [d_1,d_2]))=\lambda(A_{d_1,p}(X)\cap [d_1/N^p,d_2-d_1+d_1/N^p])$$
We prove 
$$A_{d_1,p}(X)\cap [d_1/N^p,d_2-d_1+d_1/N^p]=X\cap [d_1/N^p,d_2-d_1+d_1/N^p]$$
Let $y\in A_{d_1,p}(X)\cap [d_1/N^p,d_2-d_1+d_1/N^p]$. Then $y=A_{d_1,p}(x)$ for some\\ 
$x\in X\cap [d_1,d_2]$. But on $[d_1,d_2]$ the slope of $A_{d_1,p}$ is 1, therefore $(x,y)\in S$. 
We get $y\in X$ from the fact that $x\in X$ and $X$ is saturated. Vice-versa, let\\ 
$y\in X\cap [d_1/N^p,d_2-d_1+d_1/N^p]$. Then $y=A_{d_1,p}(x)$ where $x=A_{d_1,p}^{-1}(y)$ 
which together with $X$ being saturated insures $x\in X$ (also, the slope of $A_{d_1,p}^{-1}$ is 1, 
around $y$). In conclusion the above sets are equal. From the last two relations we obtain 
$$\lambda(X\cap [d_1,d_2])=\lambda(X\cap [d_1/N^p,d_2-d_1+d_1/N^p])$$
Taking the limit when $p$ goes to infinity we obtain (\ref{e4}). If $0<d_1<d_2<d_3<1$ are three 
consecutive N-adic numbers then from (\ref{e4}) 
$\lambda(X\cap [d_1,d_2])=\lambda(X\cap [d_2,d_3])$. For any $p\in\mathbb{N}$, covering the unit 
interval with $N^p$ consecutive N-adic rationals we obtain 
$$\lambda(X)=N^p\lambda(X\cap [d_i,d_{i+1}])=\frac{\lambda(X\cap [d_i,d_{i+1}])}{\lambda([d_i,d_{i+1}])}$$
where $d_{i+1}-d_i=1/N^p$. 
Suppose now $\lambda(X)>0$. Then there exists $x\in X$ a Lebesgue point. For any $p$ we can find 
a sequence $(\frac{k_p}{N^p})_{p>0}$ of N-adic such that $x\in\cap_{p>0}[\frac{k_p}{N^p}, \frac{k_p+1}{N^p}]$. 
Hence $N^p\lambda(X\cap [k_p/N^p, (k_p+1)/N^p])\rightarrow 1$ when $p\rightarrow\infty$. This together 
with the last equality implies $\lambda(X)=1$. In conclusion $S_N$ is ergodic. 
\par 
Let $S$ be the equivalence relation determined by the N-adic translations modulo the unit interval, i.e. 
$(x,y)\in S$ iff $|x-y|=d$ for some N-adic $d\in [0,1]$ (by remark \ref{r2}, $S$ is not included in $R_N$). Notice that if $(x,y)\in S_N$ then $f(x)=y$ with 
$f^{'}(x)=1$ and $f\in F(N)$. This implies $f(x)=x+d$ for $d$ N-adic, therefore $(x,y)\in S$. Because 
$S$ is hyperfinite we obtain $S_N$ hyperfinite: indeed, write 
the equivalence class $S[x]=\cup_n R_n[x]$ where $(R_n)_n$ is an increasing sequence of 
finite equivalence relations. Then $(S_N\cap R_n)_n$ is an increasing sequence of finite equivalence 
relations. We argue that $S_N$ is with infinite orbits: let $x\in [0,1]$ and $d<x$ a N-adic. For all 
sufficiently large $p$ we have $A_{d,p}(x)\in S_N[x]$. Now, if $A_{d,p_1}(x)=A_{d,p_2}(x)$ then, as 
$p_1$ and $p_2$ are large enough we have $x-d+d/N^{p_1}=x-d+d/N^{p_2}$ hence $p_1=p_2$. In conclusion 
the $S_N$-orbit of $x$ is infinite.
\end{proof}
Let $\phi$ be the faithful normal state determined by the scalar
product with $\varphi_0$. Recall the definition of the
centralizer\\ $M^{\phi}=\{L\in M(R) : \phi(LT)=\phi(TL),\mbox{
}\forall\mbox{ }T\in M(R)\}$. We know from \cite{Co} that for
$III_{\lambda}$ factors there exists a faithful normal state such
that the centralizer is a factor of type $II_1$. 
We are now ready to prove the main result of this section.
\begin{theorem} The von Neumann algebra $M(R_N)$ is the hyperfinite factor of type $III_{1/N}$. 
The core $M^{\phi}$ is the  hyperfinite $II_1$ factor isomorphic to $M(S_N)$. 
\end{theorem}
\begin{proof}The above proposition shows that $R_N$ is ergodic as well, therefore $M(R_N)$ is a factor. 
Suppose now that there exists a $\sigma$-finite measure $\mu$, $R_N$-invariant
such that $\mu\prec\lambda$. We take the Radon-Nikodym derivative 
$f:=\frac{\partial\mu}{\partial\lambda}$. By invariance of $\mu$ with respect to $R_N$ and with a substitution we obtain
\begin{equation}\label{e5}
f(x)=f(T(x))T^{'}(x), \mbox{ }\forall\mbox{ } T\in F(N), \mbox{a.e.} x\in [0,1]
\end{equation}
For some fixed values $a<b$, we consider the set \\
$A:=\{x\in [0,1]\mbox{ }|\mbox{ }f(x)\in [a,b]\}$. We show that $A$ is $S_N$-saturated: if 
$(x,y)\in S_N$ with $x\in A$ then there is a $T\in F(N)$ such that $T(x)=y$ and $T^{'}(x)=1$. Applying 
equation (\ref{e5}) we get $f(y)=f(x)$, therefore $y\in A$. By ergodicity $A$ has to be of Lebesgue measure 0 
or 1. Because $a$ and $b$ are arbitrary we obtain
that $f$ must be constant. This is not possible though, as the
Lebesgue measure $\lambda$ is not $R_N$-invariant. In conclusion
there is no such measure $\mu$, therefore $M(R_N)$ is a factor of type III. \\
We prove next that $M(R_N)$ is of type $III_{1/N}$. 
We use Proposition 2.2 in \cite{FM1}:
in particular it says that for $T\in F$ we have $D(T^{-1}(y),y)=dT_{*}(\lambda)/d\lambda(y)$ a.e. $y$. Thus for any Borel 
subset $A$ of $[0,1]$ we have
$$\int_{A}(T^{-1})^{'}(y)dy=\lambda(T^{-1}(A))=\int_{A}D(T^{-1}(y),y)dy$$
Therefore
\begin{equation}\label{e2}
\forall\mbox{ }T\in F\mbox{ : }D(x,T(x))=\frac{1}{T^{'}(x)},\mbox{ a.e.}x
\end{equation}
The above equation (almost) finds the range, $N^{\mathbb{Z}}$, of
the map $D:R\rightarrow\mathbb{R}_{+}$. Indeed, for any $(x,y)\in
R$ there exists an unique $T\in F(N)$ such that $T(x)=y$ (if there
are $T_1\neq T_2$ in $F(N)$ such that $T_1(x)=T_2(x)$ , then $x$ must
be $N$-adic rational, a value which we avoid by remark \ref{r1}). We 
will actually compute the asymptotic range of $D$,
$$r^{*}(D)=\{a|\mbox{ }\forall\mbox{ }V_a\mbox{ neighborhood of }a\mbox{ }\forall\mbox{ }Y\subset [0,1]\mbox{ of positive
measure }\mbox{ : }$$
$$\mbox{pr}\{(x,y)\in Y\times Y\mbox{ }|\mbox{ } D(x,y)\in V_a\}=Y\mbox{ a.e. }\}$$
Notice first that \\ $\mbox{pr}\{(x,y)\in Y\times Y\mbox{ }|\mbox{ } D(x,y)\in V_a\}=
\{x\in Y\mbox{ }|\mbox{ }\exists\mbox{ }T\in F(N)\mbox{:}D(x,T(x))\in V_a\}$\\
If $a\notin N^{\mathbb{Z}}$ then there is a neighborhood $V_a$ of
$a$ such that $V_a\cap N^{\mathbb{Z}}=\emptyset$. This combined
with $D(x,T(x))\in N^{\mathbb{Z}}$ and equation (\ref{e2}) implies
$\lambda\{x\in Y\mbox{ }|\mbox{ }\exists\mbox{ }T\in
F\mbox{:}D(x,T(x))\in V_a\}=0$, which means $a\notin r^{*}(D)$.\\
For proving the other inclusion, let $p\in\mathbb{Z}$ and $Y\subset [0,1]$ with $\lambda(Y)>0$. 
From the definition of the asymptotic range of the map $D$, suffices to show: a.e. $x\in Y$, $\exists$ 
$y\in Y$, $\exists$ $T\in F(N)$ such that $T(x)=y$ and $T^{'}(x)=N^{-p}$ (because $D(x,T(x))=1/T^{'}(x)$). 
For $S_N[Y]$, the saturation of $Y$ through $S_N$, we have $\lambda(S_N[Y])=1$. Consider the set 
$$Y_1:=\{y\mbox{ }|\mbox{ }\exists\mbox{ }x\in Y,\mbox{ }\exists T\in F(N)\mbox{ such that }T(x)=y, T^{'}(x)=N^{-p}\}$$
Then $\lambda(Y_1\setminus S_N[Y])=0$. Because $F(N)$ is countable and its elements preserve the null sets the 
following set is of measure 0, 
$$C:=\bigcup_{T\in F(N)}T^{-1}(Y_1\setminus S[Y])$$
Now, let $x\notin C$ and $x\in Y$. Choose $0<d<1$ a N-adic such that $x\in [0,d]$. Then for 
$T_1:=A_{d,p}$ we have $T_1^{'}(x)=N^{-p}$. The point $x$ not being in $C$ we obtain $T_1(x)\in S_N[Y]$, i.e. 
$\exists$ $T_2$ with $T_2^{'}(T_1(x))=1$ and $T_2(T_1(x))\in Y$; the point $y:=T_2(T_1(x))$ is the one 
we are looking for. Therefore $M(R_N)$ is a type $III_{1/N}$ factor. To check the last part of the theorem 
notice that the kernel of the Radon-Nikodym derivative $D$ equals precisely $S_N$. From here it is rather 
standard (\cite{Co}, \cite{Ta}) to conclude that the core $M^{\phi}$ is $M(S_N)$ and $M(R_N)$ is the 
crossed-product of the hyperfinite $II_{\infty}$ factor by $\mathbb{Z}$.
\end{proof}
\begin{remark}For the particular case when $N=2$, Sergey Neshveyev pointed out that we can show that
$M(R_2)$ is the hyperfinite $III_{1/2}$ factor in the following way:
\par
Consider the $ax+b$ group with $a$ of the form $2^n$,
$n\in\mathbb{Z}$ and $b\in\mathbb{R}$, dyadic, i.e. $b=k/2^m$ for
some $k,m\in\mathbb{Z}$. The multiplication is given by
$(a,b)(a',b')=(aa',ab'+b)$. This groups acts naturally on
$\mathbb{R}$ by dyadic translations and dilations by powers of
$2$. It therefore generates an equivalence relation on
$\mathbb{R}$. From remark \ref{r2}  the restriction of this
equivalence relation to $[0,1]$ is $R_2$.
\par
The crossed-product $M_{\mathbb{R}}$ of $\linfr$ with this action
decomposes as follows: first, the dyadic translations act freely
and ergodically on $\linfr$, so that the crossed-product is a hyperfinite $II_\infty$
factor. Then the dilations by $2$ induce an automorphism on this
$II_\infty$ factor that scales the semi-finite trace by $2$.
Therefore, using Connes results \cite{Co}, we get that
$M_{\mathbb{R}}$ is a hyperfinite $III_{1/2}$ factor. Now take the
projection $p$ given by the characteristic function of $[0,1]$.
The compression $pM_{\mathbb{R}}p$ is isomorphic to
$M_{\mathbb{R}}$ (since we are in a type $III$ factor); on the
other hand, it can be shown that this compression is isomorphic to
our $M(R_2)$.\\ 
Notice that the same "compression" argument cannot work for general $N$, see 
the counterexample in remark \ref{r2}.  
\end{remark}
\section{\label{tr}A treeing of $R_2$}
Let us notice that $R_N$ is treeable being a hyperfinite equivalence relation. This is a consequence of the 
general theory developed mainly by H.Dye, W.Krieger and Connes-Feldman-Weiss. In the following using the finite generation of $F$ 
we will precisely find such a treeing.\\ 
Let  $A$ and $B$ the piecewise linear homeomorphisms that generate $F$. Let us consider the following graphing:\\
$\Phi=(\varphi_i:A_i\rightarrow B_i)_{i\in\{1,2,3\}}$ where $\varphi_i$'s are defined as follows:\\
$\varphi_1:[0,1/2]\rightarrow [0,3/4]$, $\mbox{  }\varphi_1(x)=A^{-1}(x)$,\\
$\varphi_2:[1/2,3/4]\rightarrow [1/2,7/8]$, $\mbox{  }\varphi_2(x)=B^{-1}(x)$,\\
$\varphi_3:[3/4,1]\rightarrow [1/2,1]$, $\mbox{  }\varphi_3(x)=A(x)$.\\
\begin{proposition}\label{pr1} $R=R_{\Phi}$
\end{proposition}
\begin{proof} Clearly $R_{\Phi}\subset R$. Let $(x,y)\in R$ i.e. $\omega (x)=y$ for $\omega\in F$ word over the letters $A$,
$A^{-1}$, $B$, $B^{-1}$. Notice that suffices to show $(x,y)\in R_{\Phi}$ for $\omega\in\{A,B\}$ (apply induction on the
length of $\omega$). \\
 Case I:  $A(x)=y$\\
I.1: If $x\in [0,1/2]$ then $A(x)=x/2=y\in [0,1/4]\subset [0,1/2]$, hence $x=\varphi_1(y)$\\
I.2: If $x\in [1/2,3/4]$ then $A(x)=y=x-1/4\in [1/4,1/2]\subset [0,1/2]$, hence $x=\varphi_1(y)$\\
I.3: If $x\in [3/4,1]$ then $\varphi_3(x)=y$\\
For this case we conclude $(x,y)\in R_{\Phi}$.\\
 Case II: $B(x)=y$\\
II.1: If $x\in [0,1/2]$ then $x=y$\\
II.2: If $x\in [1/2,3/4]$ then $y=x/2+1/4\in [1/2,5/8]$, hence $x=\varphi_2(y)$\\
II.3: If $x\in [3/4,7/8]$ then $y=x-1/8\in [5/8,3/4]$, hence $x=\varphi_2(y)$\\
II.4: If $x\in [7/8,1]$ then $y=2x-1\in [3/4,1]$, hence $x=\varphi_3(y)$\\
From all cases we conclude $(x,y)\in R_{\Phi}$
\end{proof}
\begin{theorem}For all $\omega$ reduced words over $\Phi$, the set $\{x\in [0,1]\mbox{ }|\mbox{ }\omega (x)=x\}$ has Lebesque
measure zero, i.e. almost every orbit has a tree structure.
\end{theorem}
\begin{proof} If $\omega=\varphi_{i_1}^{\epsilon_1}\varphi_{i_2}^{\epsilon_2}...\varphi_{i_k}^{\epsilon_k}$ is a reduced word
over $\Phi$ then $i_j\in\{1,2,3\}$, $\epsilon_j\in \{-1,1\}$ and if $i_j=i_{j+1}$ then $\epsilon_j=\epsilon_{j+1}$. To show
the set of fixed points has measure zero we use induction on the length $k$. The case $k=1$ being trivial we assume for any
reduced word of length $k-1$ the measure of its fixed points is 0. Take $\omega$ of length $k$ and $x$
such that $\omega(x)=x$. We may discard the orbits of $x=1/2$ and $x=3/4$ as
these are countable sets. We distinguish three cases:\\
\\ I. $x\in [0,1/2)$ \\
We must have $i_k=1=i_1$, $\varphi_1$ being the only generator whose domain is $[0,1/2]$ and that can target points in
$[0,1/2)$. If $\epsilon_1\neq\epsilon_k$ apply the induction hypothesis for the word
$\varphi_{i_2}^{\epsilon_2}...\varphi_{i_k}^{\epsilon_k}$. If $\epsilon_1=\epsilon_k=1$ then
$\varphi_{i_2}^{\epsilon_2}...\varphi_{i_{k-1}}^{\epsilon_{k-1}}\varphi_1(x)=\varphi_1^{-1}(x)\in [0,1/2)$. As above we
obtain $i_2=1$. $\omega$ being reduced we have $\epsilon_1=\epsilon_2$ so that $\varphi_{i_3}^{\epsilon_3}...\varphi_{i_{k-1}}^{\epsilon_{k-1}}\varphi_1(x)=\varphi_1^{-2}(x)\in [0,1/2)$.
Inductively we obtain all subscripts $i_j=1$. The equation $\omega (x)=x$ becomes $\varphi_1^k(x)=x$, therefore there is
at most one solution for $\omega(x)=x$. By symmetry, the case $\epsilon_1=\epsilon_k=-1$ has a similar argument. \\
\\ II. $x\in (1/2,3/4)$\\
Suppose $i_k=1$. In order $\varphi_1^{\epsilon_k}(x)$ to make sense we must have $\epsilon_k=-1$.
Because $\omega$ is reduced and $\varphi_1^{-1}(x)\in [0,1/2)$ the only choice for the letter
$\varphi_{i_{k-1}}^{\epsilon_{k-1}}$ is $\varphi_1^{-1}$. Continuing this procedure we would make all letters of $\omega$
equal to $\varphi_1^{-1}$, i.e. $x$ is a fixed point of $\varphi_1^{k}$. Same conclusion holds if $i_1=1$. Suppose now
$i_1$, $i_k\in\{2,3\}$. We distinguish the following subcases:\\
II.1. $\omega (x)=\varphi_2^{\epsilon_1}\overline{\omega}\varphi_2^{\epsilon_k}(x)=x$. If $\epsilon_1\neq\epsilon_k$ then the induction
hypothesis will end the proof. By symmetry suffices to check only the case $\varphi_2^{-1}\overline{\omega}\varphi_2^{-1}(x)=x$.
We claim that $\overline{\omega}$ has the following $\Phi$-writting:
$\overline\omega=\varphi_3^{-p_1}\varphi_2^{-q_1}...\varphi_3^{-p_l}\varphi_2^{-q_l}$ where $p_j\geq 0$, $q_j\geq 0$ are integers.
Because of the way we choose the domains the following statements are true (the "reading" of $\omega$ is done from right to
left, i.e. letter $x$ is after letter $y$ in $xy$): \\
$\bullet$ There can be no $\varphi_1^{\epsilon}$ occurence in $\overline\omega$: inded,
a $\varphi_1$ occurence will force all letters to the right
of $\varphi_1$ be equal to $\varphi_1$. This in not allowed as the right-end letter takes on $x\in [1/2,3/4]$. A $\varphi_1^{-1}$
occurence is not allowed otherwise all letters to the left of it would be equal to $\varphi_1^{-1}$, including the left-end.
In this case $\omega(x)=x$ would be sent in $[0,1/2]$.\\
$\bullet$ A $\varphi_3$ occurence immediately after $\varphi_2^{-1}$ is not possible.\\
$\bullet$ After a $\varphi_2^{-1}$ occurence only a $\varphi_3^{-1}$ or $\varphi_2^{-1}$ occurence is allowed.\\
$\bullet$ A $\varphi_2$ occurence immediately after $\varphi_3^{-1}$ is not possible.\\
$\bullet$ After a $\varphi_3^{-1}$ occurence only a $\varphi_3^{-1}$ or $\varphi_2^{-1}$ occurence is allowed.\\
All of the above prove the claim. We show that the equation $\omega (x)=x$ has at most one solution: $\varphi_2^{-1}$ takes
$[1/2,3/4]$ into $[1/2,3/4]$ and $\varphi_2^{-1}(x)=x/2+1/4$ so that with each iteration the slope will decrease by a factor
of $2$; we apply a $1/2$ slope at least once at the right-end of $\omega$ when computing $\varphi_2^{-1}(x)$ (it may be that at some
step in the composition the trajectory exits $[1/2,3/4]$ and $\varphi_2^{-1}$ takes on slope $=1$ but the slope has already
been "damaged" at the begining); the slope is decreased further by $\varphi_3^{-1}=(x+1)/2$. Now the equation $\varphi_2^{-1}\overline{\omega}\varphi_2^{-1}(x)=x$
can be written $ax+b=x$ for some $a<1$.\\
II.2. $\omega (x)=\varphi_3^{\epsilon_1}\overline{\omega}\varphi_3^{\epsilon_k}(x)=x$. Again by the induction hypothesis
suffices to argue only for the case $\varphi_3^{-1}\overline{\omega}\varphi_3^{-1}(x)=x$: this is easy as $\varphi_3^{-1}$
targets $[3/4,1]$ but $x\in (1/2,3/4)$. \\
II.3. $\omega (x)=\varphi_2^{\epsilon_1}\overline{\omega}\varphi_3^{\epsilon_k}(x)=x$. Because $x<3/4\mbox{ : }$
$\epsilon_k=-1$. A similar analysis of occurences and slopes$<1$ leads to an equation with one solution at most.\\
II.4.  $\omega (x)=\varphi_3^{\epsilon_1}\overline{\omega}\varphi_2^{\epsilon_k}(x)=x$. This case is symmetric to II.3.\\
\\ III. $x\in (3/4,1]$\\
Again we discard $\varphi_1^{\epsilon}$'s occurences in $\omega$: a $\varphi_1$ occurence will force all letters to the right
of $\varphi_1$ be equal to $\varphi_1$. This in not allowed as the right-end letter takes on $x\in (3/4,1]$. A $\varphi_1^{-1}$
occurence is not allowed otherwise all letters to the left of it would be equal to $\varphi_1^{-1}$, including the left-end.
In this case $\omega(x)$ would be sent in $[0,1/2]$. We list now all possibilities for the first and last letter of $\omega$: \\
$\varphi_{i_1}^{\epsilon_1}\in\{\varphi_2,\mbox{ }\varphi_3^{\epsilon}\}$,
$\varphi_{i_k}^{\epsilon_k}\in\{\varphi_2^{-1},\mbox{ }\varphi_3^{\epsilon}\}$ where $\epsilon\in\{-1,1\}$.\\ The cases
$\omega\in\{\varphi_3^{-1}\overline\omega\varphi_3,\mbox{ }\varphi_2\overline\omega\varphi_2^{-1},\mbox{ }\varphi_3\overline\omega\varphi_3^{-1}\}$
can be dealt with by the induction hypothesis. All the other remaining cases can be dealt with by the same analysis of
occurences in II.1 : e.g. if
$\omega= \varphi_3^{-1}\overline\omega\varphi_3^{-1}$ then the first letter (from the right) of $\overline\omega$ is either
$\varphi_3^{-1}$ or $\varphi_2^{-1}$ etc; in the end $\overline\omega$ becomes a word written with iterates
of $\varphi_3^{-1}$ and/or $\varphi_2^{-1}$. Because $\varphi_3^{-1}$ has slope $1/2$ and $\varphi_2^{-1}$ has slope $1$ or
$1/2$ we conclude that the equation $\omega(x)=x$ is equivalent to $ax+b=x$ with $a<1$. \\
With the analysis of I, II and III we complete the $k^{\mbox{th}}$ step of induction, thus proving the theorem.
\end{proof}

\begin{remark} Using the infinite presentation of the Thompson group $F$ it can be shown $C(F)=1$. Using 
Gaboriau's results we will describe how to compute this cost, but for the general version $F(N)$. 
Still, the question is whether the cost of the normal subgroup $[F,F]$ is 1 or $>1$(in this case $F$ 
would be non-amenable); we believe it should be 1, even though we do not know if the following 
procedure can be carried-out for $[F,F]$ instead. 
\end{remark}
\par 
The following properties are easy to work-out:\\
i) any non-trivial element of $F(N)$ is of infinite order;\\
ii) $x_Nx_1^{-1}$ commutes with any $x_j\in F(N)$, where $j>N$.\\
\begin{proposition} $C(F(N))=1$
\end{proposition}
\begin{proof}The idea of the proof is similar to the case $N=2$ which is done in \cite{Gab}. We first 
show that the group $\Gamma$ generated by $\gamma:=x_Nx_1^{-1}$ and $x_i$, 
$i>N$ has fixed price $=1$. Let $\Pi:\Gamma\rightarrow\mbox{Aut}(X,\nu)$ be a
free action that generates a (SP1) equivalence relation $R_{\Pi}$ of $\Gamma$. We prove $C(R_{\Pi})=1$: suffices to show for every
$\delta>0$, $C(R_{\Pi})\leq 1+\delta$. Because $\gamma$ is of infinite order we can find a sequence $A_n$ of Borel subsets of
$X$ such that $\nu(A_n)<\delta/2^n$ and $A_n\cap R_{\gamma}[x]\neq\emptyset$ a.e. $x\in X$ (see \cite{Gab}). Using ii) 
above, it is a routine to show that for the following graphing $\Phi$ we have $R_{\Phi}=R_{\Pi}$:\\
$\Phi:=\{\Pi(\gamma):X\rightarrow X,\mbox{ }\Pi(x_i)_{|A_i}\mbox{ }i>N\}$.  
Next, take $\Gamma_1$ the subgroup generated by $\Gamma$ and $x_1$. It is easy to see that the set 
$x_1\Gamma x_1^{-1}\cap\Gamma$ is infinite (it contains all $x_j$ with $j>N$. Inductively, in $N$ 
steps we obtain an increasing sequence of subgroups whose union equals $F(N)$. We apply now 
Critere 3 in \cite{Gab} to conclude that the cost of $F(N)$ is 1.
\end{proof} 
The reason all of the above does not work for the subgroup $[F,F]$
is that we do not know the generators of $[F,F]$. We can still
start with the element $\gamma$ and then gradually add elements of
$[F,F]$, the idea being to enter the hypotheses of $Critere$ 3:
however we did not find a way of adding such that to exhaust
$[F,F]$.
\begin{acknowledgements}We would like to thank prof. Florin Radulescu, Sergey Neshveyev and Ionut Chifan for useful comments and suggestions about the
subject.  We thank the referee for many simplifications of the proofs and for pointing out the right 
bibliography. The second author also thanks prof. Florin Radulescu for his careful guidance and support 
along the years.
\end{acknowledgements}


\begin{thebibliography}{abcdefg}
\bibitem[Br]{Br}
K.S.Brown, {\em {F}initeness {P}roperties of {G}roups}, J. Pure Appl. Algebra, 44 (1986), 45-75
\bibitem[Can]{Can}
J.W.Cannon, W.J.Floyd, and W.R.Parry, {\em {I}ntroductory {N}otes on {R}ichard {T}homson's {G}roups}, L'Enseignement
Mathematique, t.42 (1996), p.215-256
\bibitem[Co]{Co}
A.Connes, {\em {U}ne classification des facteurs de type III}, Ann. scient.Ec.Norm.Sup. $4^{e}$ serie, t.6, 1973, p.133-252
\bibitem[FMI]{FM1}
J.Feldman, C.Moore, {\em {E}rgodic {E}quivalence {R}elations, {C}ohomology, and {V}on {N}eumann {A}lgebras.I},
Transactions of the AMS, dec. 1977, vol.234, issue 2, p.289-324
\bibitem[FMII]{FM}
J.Feldman, C.Moore, {\em {E}rgodic {E}quivalence {R}elations, {C}ohomology, and {V}on {N}eumann {A}lgebras.II},
Transactions of the AMS, dec. 1977, vol.234, issue 2, p.325-359
\bibitem[Gab]{Gab}
D.Gaboriau, {\em {C}o$\hat{u}$t des relations d'equivalence et des groupes}, Invent.Math.139, 41-98 (2000)
\bibitem[Kr]{Kr}
W.Krieger, {\em {O}n {C}onstructing {N}on$-^{*}${I}somorphic {H}yperfinite {F}actors of {T}ype III},
Journal of Functional Analysis 6, 97-109 (1970)
\bibitem[OW]{OW}
D.Ornstein, B.Weiss, {\em {E}rgodic theory of amenable group actions I. The Rohlin lemma}, Bull.A.M.S. 2, 161-164, 1980
\bibitem[Ta]{Ta}
M.Takesaki, {\em {S}tructure of {F}actors and {A}utomorphism {G}roups}, CBMS Regional Conference Series in Mathematics, 51. AMS,1983 
\end{thebibliography}
\end{document}